\documentclass[oneside,english]{article}
\usepackage[T1]{fontenc}
\usepackage[latin1]{inputenc}
\usepackage{amssymb,amsthm,amsmath}

 \theoremstyle{plain}
 \newtheorem{thm}{Theorem}[section]
 \numberwithin{equation}{section}
 \numberwithin{figure}{section}
 \theoremstyle{plain}
 \theoremstyle{remark}
 \newtheorem{rem}[thm]{Remark}
 \theoremstyle{plain}
 \newtheorem{fact}[thm]{Fact}
 \theoremstyle{plain}
 \newtheorem{lem}[thm]{Lemma}
 \theoremstyle{remark}
 \newtheorem*{acknowledgement*}{Acknowledgement}

\usepackage{babel}
\begin{document}

\title{Occupation time fluctuations of Poisson and equilibrium finite variance
branching systems}

\author{Piotr Mi\l{}o\'s\\
Institute of Mathematics\\
Polish Academy of Sciences\\Warsaw}

\maketitle
\begin{abstract}

Functional limit theorems are presented for the rescaled occupation time
fluctuation process of a critical finite variance branching particle system in $\mathbb{R}^{d}$
with symmetric $\alpha$-stable motion starting off from either a standard Poisson
random field or from the equilibrium distribution for intermediate dimensions
$\alpha<d<2\alpha$. The limit processes are determined by sub-fractional and fractional
Brownian motions, respectively.
\end{abstract}

AMS subject classification: primary 60F17, 60G20, secondary 60G15\\

Key words: Functional central limit theorem; Occupation time fluctuations; Branching
particles systems; Fractional Brownian motion; Sub-fractional Brownian
motion; equilibrium distribution.

\section{Introduction}

Consider a system of particles in $\mathbb{R}^{d}$ starting off at
time $t=0$ from a certain distribution (a standard Poisson 
and equilibrium fields are investigated in this paper). They evolve independently,
moving according to a symmetric $\alpha$-stable L\'evy process and
undergoing finite variance branching at rate $V$ ($V>0$). We obtain functional
limit theorems for the rescaled occupation time fluctuations of this system
when $\alpha<d<2\alpha$. This is an extension of \cite[Theorem 2]{BGT2}
where the starting distribution is a Poisson field and the branching
law is critical and binary.

\subsection{Branching law\label{sub:Branching-law}}

~\\
In the \cite{BGT1,BGT2,BGT3} the law of branching is critical and
binary. In this paper an extended model is investigated. The particles
branch according to the law given by a moment generating function
$F$. $F$ fulfills two requirements: 

\begin{enumerate}
\item $F'(1)=1$, which means that the law is critical (the expected number
of particles spawning from one particle is $1$), 
\item $F''(1)<+\infty$, which states that the second moment exists. 
\end{enumerate}
(Note here that the branching law in \cite{BGT2} is given by $F(s)=\frac{1}{2}\left(1+s^{2}\right)$
and obviously fulfills the two requirements.) Although constraints
imposed on $F$ are not very restrictive and quite natural (so that
the class of the branching laws satisfying them is broad) still there
remain other interesting cases to be investigated. One of them is
the class of branching laws in the domain of attraction of the $\left(1+\beta\right)-$stable
law (i.e., the moment generating function is $F\left(s\right)=s+\frac{1}{2}\left(1+s\right)^{1+\beta}$),
the case studied in \cite{BGT4,BGT5}. A remarkable feature of the
latter case is that the limit processes are stable ones and not Gaussian
as it occurs in the finite variance case.

\subsection{Equilibrium distribution}
~\\
Another concept naturally related to particle systems is an equilibrium
distribution. It has been shown that in certain circumstances the
system converges to the equilibrium distribution \cite{GW2}. It is
both an interesting and important question whether the theorems shown
by Bojdecki et al still hold in the case when the equilibrium state
is taken as the initial condition. A conjecture in \cite{BGT1}
states that the temporal structure of the limit is given by fractional
Brownian motion. It is of interest to notice that the limit is different
from the one in the case of the system starting off from the Poisson
field (where temporal structure is sub-fractional Brownian motion).
We study behavior of the system for a branching law given by $F$.
But there is still broad area for further studies. No attempt has
been made to develop more general theory concerning systems with a
general starting distribution (or a large class of distributions).

\subsection{General concepts and notation\label{sub:General-concepts-and}}

~\\
Let us denote $N_{t}^{Poiss}$ and $N_{t}^{eq}$, the empirical processes
for the system starting off from the Poisson field with Lebesgue intensity measure and the equilibrium
respectively. For a measurable set $A\subset\mathbb{R}^{d}$, $N_{t}^{Poiss}\left(A\right)$, $N_{t}^{eq}\left(A\right)$, respectively
are the numbers of particles of the system in set
$A$ at time $t$. Note that they are measure-valued processes but we will
consider them as processes with values in $\mathcal{S}'$ (the space
of tempered distributions) because this space has good analytical
properties.

The equilibrium distribution is defined by \[
\lim_{t\rightarrow+\infty}N_{t}^{Poiss}=N_{eq},\]
where the limit is understood in weak sense.
The Laplace functional of the equilibrium distribution is given by
\begin{equation}
	\mathbb{E}\exp\left\{ -\left\langle N_{eq},\varphi \right\rangle \right\} =\exp\left\{ \left\langle \lambda,e^{-\varphi}-1\right\rangle +V\int_{0}^{\infty}\left\langle \lambda,H\left(j\left(\cdot,s\right)\right)\right\rangle ds\right\} ,\label{eq: laplace_equilibrum}\end{equation}
where 
\begin{equation}
	j\left(x,l\right):=\mathbb{E}\exp\left(-\left\langle N_{l}^{x},\varphi \right\rangle \right)
	\label{def: h}
\end{equation}
 $H(s)=F\left(s\right)-s$, $\varphi:\mathbb{R}^d\rightarrow\mathbb{R}_{+}$, $\varphi \in \mathcal{L}^1(\mathbb{R}^d)\cap C(\mathcal{R}^d)$ and $j$ satisfies the integral equation\[
j\left(x,l\right)=\mathcal{T}_{l}e^{-\varphi}\left(x\right)+V\int_{0}^{l}\mathcal{T}_{l-s}H\left(j\left(\cdot,s\right)\right)\left(x\right)ds,\]
 This equations can be obtained in the same way as \cite[(2.4)]{GW2}. Note that in \cite{GW2} function $\varphi$ is continuous with compact support. We approximate $\varphi\in \mathcal{L}^1$ using functions $\varphi_n$ with compact support $\varphi_n\nearrow \varphi$. Using Lebesgue's monotone convergence theorem it is easy to obtain the above equations for $\varphi$ ($H$ is decreasing because of the criticality of the branching law).

For an empirical process $N_{t}$ the rescaled occupation time fluctuation
process is defined by\begin{equation}
X_{T}\left(t\right)=\frac{1}{F_{T}}\int_{0}^{Tt}\left(N_{s}-\mathbb{E}N_{s}\right)ds,\: t\geq0,\label{eq occupation}\end{equation}
where $T>0$ and $F_{T}$ is a suitable
norming. We are interested in the weak functional limit of $X_T$ when time is accelerated (i.e., $T$ tends to $\infty$).  

The $\alpha$-stable process starting from $x$ will be denoted by
$\eta_{t}^{x}$ its semigroup by $\mathcal{T}_{t}$ and its infinitesimal
operator by $\Delta_{\alpha}$. The Fourier transform of $\mathcal{T}_{t}$
is \begin{equation}
\widehat{\mathcal{T}_{t}}\varphi\left(z\right)=e^{-t\left|z\right|^{\alpha}}\widehat{\varphi}\left(z\right).\label{eq: Fourier T_t}\end{equation}
For brevity let us denote \begin{equation}
K=\frac{V\Gamma\left(2-h\right)}{2^{d-1}\pi^{d/2}\alpha\Gamma\left(d/2\right)h\left(h-1\right)},\label{def: K}\end{equation}
where 
\begin{equation}
h=3-d/\alpha
\end{equation}
(in this paper we always assume that $\alpha<d<2\alpha$ so $h>1$) and
\begin{equation}
	M=F''\left(1\right).\label{eq: M}
\end{equation}
We will now introduce two centered Gaussian processes. One of them
is sub-fractional Brownian motion with parameter $h$ with the
covariance function $C_{h}$\begin{equation}
C_{h}\left(s,t\right)=s^{h}+t^{h}-\frac{1}{2}\left[\left(s+t\right)^{h}+\left|s-t\right|^{h}\right]\label{def: C - subfractional}\end{equation}
and the second one is fractional Brownian motion with parameter
$h$ and the covariance function $c_{h}$ \begin{equation}
c_{h}\left(s,t\right)=\frac{1}{2}\left(s^{h}+t^{h}-\left|s-t\right|^{h}\right).\label{def: c - fractional}\end{equation}

\subsection{Space-time method}

~\\
The space-time method is a very convenient technique for investigating
the weak convergence in the $C\left(\left[0,\tau\right],\mathcal{S}'\left(\mathbb{R}^{d}\right)\right)$
space. It was developed by Bojdecki et al and can be found in \cite{BGR}.
If $X=\left(X\left(t\right)\right)_{t\in\left[0,\tau\right]}$ is
a continuous $\mathcal{S}'\left(\mathbb{R}^{d}\right)$-valued process
we define a random element $\tilde{X}$ of $\mathcal{S}'\left(\mathbb{R}^{d+1}\right)$
by
\begin{equation}
	\left\langle \tilde{X},\Phi\right\rangle =\int_{0}^{\tau}\left\langle X\left(t\right),\Phi\left(\cdot,t\right)\right\rangle dt,\label{eq:space-time}
\end{equation}
where $\Phi\in\mathcal{S}\left(\mathbb{R}^{d+1}\right)$. In order
to prove that $X_{T}$ converges weakly to $X$ in $C\left(\left[0,\tau\right],\mathcal{S}'\left(\mathbb{R}^{d}\right)\right)$
it suffices to show that \[
\left\langle \tilde{X}_{T},\Phi\right\rangle \Rightarrow\left\langle \tilde{X},\Phi\right\rangle ,\,\forall_{\Phi\in\mathcal{S}\left(\mathbb{R}^{d+1}\right)}\]
 and that the family $X_{T}$ is tight.

\section{Convergence theorems}

We will present two theorems. In the first of them (which is a direct
extension of \cite[Theorem 2.2]{BGT2}) we study the occupation time
fluctuation process for the branching system starting off from the
Poisson field with Lebesgue intensity measure (denoted by $\lambda$) with the branching law given by a moment generating function as described in Section \ref{sub:Branching-law}. The result is very similar
to the one obtained in \cite[Theorem 2.2]{BGT2} - namely, the limit
process is the same up to constants. 

\begin{thm}
\label{thm:For-the-branching}Assume that $\alpha<d<2\alpha$ and
let $X_{T}$ be the occupation time fluctuation process defined by (\ref{eq occupation}) for the branching
system $N^{Poiss}$, and $F_{T}=T^{\left(3-\frac{d}{\alpha}\right)/2}$.
Then $X_{T}\Rightarrow X$ in $C\left(\left[0,\tau\right],\mathcal{S}'\left(\mathbb{R}^{d}\right)\right)$
as $T\rightarrow+\infty$ for any $\tau>0$, where $\left(X\left(t\right)\right)_{t\geq0}$
is a centered $\mathcal{S}'$-valued, Gaussian process with covariance
function: \begin{equation}
Cov\left(\left\langle X\left(s\right),\varphi\right\rangle ,\left\langle X\left(t\right),\psi\right\rangle \right)=KM\left\langle \lambda,\varphi\right\rangle \left\langle \lambda,\psi\right\rangle C_{h}\left(s,t\right),\label{def:kowariancja_tw1}\end{equation}
where$\ \varphi,\psi\in\mathcal{S}\left(\mathbb{R}^{d}\right)$.
\end{thm}
The second theorem concerns the case where the system starts from the equilibrium distribution. As
it was mentioned hereinabove the theorem is interesting because the
limit has a different time structure from the one in \cite[Theorem 2.2]{BGT2}
and Theorem \ref{thm:For-the-branching}.

\begin{thm}
\label{thm:equlibiurim}Assume that $\alpha<d<2\alpha$ and let $X_{T}$
be the occupation time fluctuation process defined by (\ref{eq occupation}) for the branching system
$N^{eq}$, and $F_{T}=T^{\left(3-\frac{d}{\alpha}\right)/2}$.
Then $X_{T}\Rightarrow X$ in $C\left(\left[0,\tau\right],\mathcal{S}'\left(\mathbb{R}^{d}\right)\right)$
as $T\rightarrow+\infty$ for any $\tau>0$, where $\left(X\left(t\right)\right)_{t\geq0}$
is a centered Gaussian process with the covariance function \begin{equation}
Cov\left(\left\langle X\left(s\right),\varphi\right\rangle ,\left\langle X\left(t\right),\psi\right\rangle \right)=KM\left\langle \lambda,\varphi\right\rangle \left\langle \lambda,\psi\right\rangle c_{h}\left(s,t\right),\label{def:kowariancja_tw2}\end{equation}
where $\varphi,\psi\in\mathcal{S}\left(\mathbb{R}^{d}\right)$.
\end{thm}
\begin{rem}
The limit processes above can be represented as follows: \\
For Theorem \ref{thm:For-the-branching}\[
X=\left(MK\right)^{1/2}\lambda\beta^{h}\]
 and for Theorem \ref{thm:equlibiurim}\[
X=\left(MK\right)^{1/2}\lambda\xi^{h},\]
where $\beta^{h}$ and $\xi^{h}$ are respectively sub-fractional and
fractional Gaussian processes defined in Section \ref{sub:General-concepts-and}.
In both cases the limit process $X$ has a trivial spatial structure
(Lebesgue measure), whereas the time structure is complicated, with
long range dependence. 
\end{rem}

\begin{rem}
The occupation time fluctuation processes of particle systems form
an area that receives a lot of research attention. We would like to
mention some other related work. Firstly the case of non-branching
systems has been studied in \cite[Theorem 2.1]{BGT2}. The result
is analogous, both to Theorem \ref{thm:For-the-branching} and \ref{thm:equlibiurim}
because the Poisson field is the equilibrium distribution for the
system. The limit process is essentially the same as in Theorem \ref{thm:equlibiurim}.
For the critical $d=2\alpha$ and large dimensions $d>2\alpha$, there
is no long range dependence and the results can be found in \cite{BGT3}.
In \cite{BZ} the fluctuations of the occupation time of the origin
are studied for a critical binary branching random walks on the $d$-dimensional
lattice, $d\geq3$, including also the equilibrium case. The convergence results are analogous to those
in \cite{BGT2,BGT3} and in this paper but the proofs are substantially different. A similar model with $\alpha=2$ was investigated in \cite{DW} (ie. with particles moving according to Brownian motion).
\end{rem}

\section{Proofs}

The main idea used in both of the proofs is to study the Laplace functional
of a process given by the space-time method. The Fourier transform
is used for this purpose. This is similar to the method in \cite{BGT2}.
In the case of Theorem \ref{thm:For-the-branching} the proof follows
the same principle as \cite[Theorem 2.2]{BGT2}. The moment generating
function can be represented using Taylor's expansion and two following
statements need to be proved. Firstly, one has to check that the method
used in \cite{BGT2} can still be applied. Secondly, it needs to be
shown that terms of order higher then $2$ play no role in the limit.
The proof of Theorem \ref{thm:equlibiurim} requires more work. The
Laplace formula contains a function that is a solution of a differential
equation. This makes the computations more cumbersome. Some expressions
in this proof had to be examined more carefully than in Theorem \ref{thm:For-the-branching}.
It should be noted that Theorem \ref{thm:equlibiurim} covers all
branching laws described in Section \ref{sub:Branching-law}.\\
Now we introduce some notation and facts used further on.

For a generating function $F$ we define \begin{equation}
G\left(s\right)=F\left(1-s\right)-1+s.\label{def: G}\end{equation}
 The following fact describes basic properties of $G$ which are straightforward
consequences of the properties of $F$.

\begin{fact}
\label{fact: properties of G}
\begin{enumerate}
\item $G\left(0\right)=F\left(1\right)-1=0$,
\item $G'\left(0\right)=-F'\left(1\right)+1=0$ since $F'\left(1\right)=1$,
\item $G''\left(0\right)=F''\left(1\right)<+\infty$,
\item $G\left(v\right)=\frac{M}{2}v^{2}+g\left(v\right)v^{2}$ where $M$
is defined by (\ref{eq: M}) and $\lim_{v\rightarrow0}g\left(v\right)=0$.
\end{enumerate}
\end{fact}
The next simple fact will be useful in proving some inequalities

\begin{fact}
\label{G>=3D0}$G\left(v\right)\geq0$ for $v\in\left[0,1\right]$. 
\end{fact}
\begin{proof}
$F''\left(1-v\right)\geq0$ which is an obvious consequence of the
fact that all of the coefficients in the expansion of $F''$ are non-negative
and $1-v\in\left[0,1\right]$. $G''\left(v\right)=F''\left(1-v\right)\geq0$.
We also know that $G'\left(0\right)=0$ so $G'\left(v\right)\geq0$
for $v\in\left[0,1\right]$. The proof is complete since $G\left(0\right)=0$
and $G$ is non-decreasing.
\end{proof}
The existence of the second moment of the moment generating function
$F$ implies also that $G$ is comparable with function $v^{2}$.

\begin{fact}
\label{G(v)/v2}We have \[
\sup_{v\in\left[0,1\right]}\frac{G\left(v\right)}{v^{2}}<+\infty\]
 
\end{fact}
\begin{proof}
Since both $G\left(v\right)$ and $v^{2}$ are continuous we only
have to check that the limit of the quotient at $v=0$ is finite.
This becomes obvious when we recall Taylor's expansion of $G\left(v\right)$
from Fact \ref{fact: properties of G}, property 4.
\end{proof}
Let us now introduce some notation used throughout the rest of the
paper. $\Phi$ will denote a positive function from $\mathcal{S}\left(\mathbb{R}^{d+1}\right)$.
\cite[Lemma in Section 3.2]{BGT2} explains why without loss of generality
it can be assumed $\Phi\geq0$. We denote\[
\Psi\left(x,s\right)=\int_{s}^{1}\Phi\left(x,t\right)dt,\]
\[
\Psi_{T}\left(x,s\right)=\frac{1}{F_{T}}\Psi\left(x,\frac{s}{T}\right).\]
To make computations less cumbersome we will sometimes assume that
$\Phi$ is of the form $\Phi\left(x,t\right)=\varphi\left(x\right)\psi\left(t\right)$
for $\varphi\in\mathcal{S}\left(\mathbb{R}^{d}\right),\,\psi\in\mathcal{S}\left(\mathbb{R}\right)$
and hence 
\begin{equation}
	\Psi_{T}\left(x,t\right)=\varphi_{T}\left(x\right)\chi_{T}\left(t\right),
	\label{eq:uproszczenie}
\end{equation}
 where $\varphi_{T}\left(x\right)=\frac{1}{F_{T}}\varphi\left(x\right)$,  $\chi\left(t\right)=\int_{t}^{1}\psi\left(s\right)ds$, $\chi_{T}=\chi\left(\frac{t}{T}\right)$.
Notice that $\varphi\geq0,\chi\geq0$ as $\Phi\geq0$.\\
Let us introduce now an important function which will appear as a
part of the Laplace functional of the occupation time fluctuation
processes \[
v_{\Psi}\left(x,r,t\right)=1-\mathbb{E}\exp\left\{ -\int_{0}^{t}\left\langle N_{s}^{x},\Psi\left(\cdot,r+s\right)\right\rangle ds\right\} ,\]
where $N_{s}^{x}$ denotes the empirical measure of the particle system
with the initial condition $N_{0}^{x}=\delta_{x}$. Let us note here that
due to the fact that $\Psi\geq0$ we have $v_{\Psi}\in\left[0,1\right]$.
We also write\begin{equation}
n_{\Psi}\left(x,r,t\right)=\int_{0}^{t}\mathcal{T}_{t-s}\Psi\left(\cdot,r+t-s\right)\left(x\right)ds.\label{def: n}\end{equation}
For simplicity of notation, we write \begin{equation}
v_{T}\left(x,r,t\right)=v_{\Psi_{T}}\left(x,r,t\right),\label{eq:v_T notacja}\end{equation}
 \begin{equation}
n_{T}\left(x,r,t\right)=n_{\Psi_{T}}\left(x,r,t\right),\label{eq:n_T notacja}\end{equation}
 \begin{equation}
v_{T}\left(x\right)=v_{T}\left(x,0,T\right),\label{eq: v_T_x notacja}\end{equation}
 \begin{equation}
n_{T}\left(x\right)=n_{T}\left(x,0,T\right)\label{eq: n_T_x notacja}\end{equation}
 when no confusion can arise.

Now we obtain an integral equation for $v$ which will play a
crucial role in the next proofs. Note that similar computations can be found also in \cite{GW1}.

\begin{lem}
$v_{\Psi}$ satisfies the equation\begin{equation}
v_{\Psi}\left(x,r,t\right)=\int_{0}^{t}\mathcal{T}_{t-s}\left[\Psi\left(\cdot,r+t-s\right)\left(1-v_{\Psi}\left(\cdot,r+t-s,s\right)\right)-VG\left(v_{\Psi}\left(x,r+t-s,s\right)\right)\right]\left(x\right)ds.\label{calkowe}\end{equation}

\end{lem}
\begin{proof}
Firstly let us investigate

\[
w\left(x,r,t\right)\equiv w_{\Psi}\left(x,r,t\right)=\mathbb{E}\exp\left(-\int_{0}^{t}\left\langle N_{s}^{x},\Psi(\cdot,r+s)\right\rangle ds\right)=1-v_{\Psi}\left(x,r,t\right),\]
We assume $\Psi\geq0$ hence we have $w\left(x,r,t\right)\in\left[0,1\right]$.
By conditioning on the time of the first branching we obtain the following
equation \begin{align*}
w\left(x,r,t\right)= & e^{-Vt}\mathbb{E}\left(-\int_{0}^{t}\Psi(\eta_{s}^{x},r+s)ds\right)\\
 & +V\int_{0}^{t}e^{-Vs}\mathbb{E}\exp\left(-\int_{0}^{s}\Psi(\eta_{u}^{x},r+u)du\right)F\left(w\left(\eta_{s}^{x},r+s,t-s\right)\right),\end{align*}
where $t\geq0,\, r\geq0$.\\
Using Feynman-Kac formula one can obtain the following equation for $w$ (for details see \cite[(3.13)-(3.17)]{BGT2})
\[
\left\{ \begin{array}{l}
\frac{\partial}{\partial t}w\left(x,r,t\right)=\left(\Delta_{\alpha}+\frac{\partial}{\partial r}-\Psi\left(x,r\right)\right)w\left(x,r,t\right)+V\left[F\left(w\left(x,r,t\right)\right)-w\left(x,r,t\right)\right],\\
w\left(x,r,0\right)=1.\end{array}\right.\]
$v\left(x,r,t\right)=v_{\Psi}\left(x,r,t\right)=1-w_{\Psi}\left(x,r,t\right)$
so $v$ satisfies the equation\[
\left\{ \begin{array}{l}
\frac{\partial}{\partial t}v\left(x,r,t\right)=\left(\Delta_{\alpha}+\frac{\partial}{\partial r}\right)v\left(x,r,t\right)+\Psi\left(x,r\right)\left(1-v\left(x,r,t\right)\right)-VG\left(v\left(x,r,t\right)\right),\\
v\left(x,r,0\right)=0.\end{array}\right.\]
 Its integral version is \ref{calkowe} (note that in \cite{BGT2} $G\left(t\right)=\frac{1}{2}t^{2}$).
\[
v\left(x,r,t\right)=\int_{0}^{t}\mathcal{T}_{t-s}\left[\Psi\left(\cdot,r+t-s\right)\left(1-v\left(\cdot,r+t-s,s\right)\right)-VG\left(v\left(x,r+t-s,t\right)\right)\right]\left(x\right)ds.\]

\end{proof}
\begin{fact}
~\begin{equation}
v_{\Psi}\left(x,r,t\right)\leq n_{\Psi}\left(x,r,t\right)\label{ineq: v<n}\end{equation}

\end{fact}
\begin{proof}
This is a direct consequence of the equation (\ref{calkowe}), the
fact that $1\geq v\geq0$ and Fact \ref{G>=3D0}.
\end{proof}
\begin{fact}
For the system $N_{t}^{Poiss}$ the covariance function is given by\begin{align}
Cov\left(\left\langle N_{u}^{Poiss},\varphi\right\rangle ,\left\langle N_{v}^{Poiss},\psi\right\rangle \right) & =\left\langle \lambda,\varphi\mathcal{T}_{v-u}\psi\right\rangle \label{eq:kowariancja Poisson}\\
 & F''\left(1\right)\cdot V\int_{0}^{u}\left\langle \lambda,\varphi\mathcal{T}_{u+v-2r}\psi\right\rangle dr,\,\,\, u\leq v,\nonumber \end{align}
where$\ \varphi,\psi\in\mathcal{S}\left(\mathbb{R}^{d}\right)$.
\end{fact}
The proof of the fact follows from a simple computation which can
be carried on using \cite[formula (3.14)]{GR}, therefore we omit it.

\subsection{Proof of theorem \ref{thm:For-the-branching}}
\subsubsection{Tightness}
The first step required to establish the weak convergence is to prove
tightness of $X_{T}$. By the Mitoma theorem \cite[Mitoma 1983]{M}
it is sufficient to show tightness of the real processes $\left\langle X_{T},\phi\right\rangle $
for all $\phi\in\mathcal{S}\left(\mathbb{R}^{d}\right)$. This can
be done using a criterion \cite[Theorem 12.3]{B}. Detailed examination
of the proof in \cite{BGT2} reveals that only the covariance function
of the $N_{t}^{Poiss}$ is needed \cite[Section 3.1]{BGT2}. One can
see that the covariance function (\ref{eq:kowariancja Poisson}) is
essentially the same as for the binary branching. Hence the proof
from \cite{BGT2} still holds for the new family of processes.

\subsubsection{The Laplace functional}

~\\
The second step uses the space-time method. According to (\ref{eq:space-time})
we define $\tilde{X}_{T}$ (from now on $\tau=1$). To establish the
convergence we use Laplace functional. By the Poisson initial condition
we have (this equation is the same as \cite[(3.10)]{BGT2}
\begin{equation}
\mathbb{E}\exp\left\{ -\left\langle \tilde{X}_{T},\Phi\right\rangle \right\} =\exp\left\{ \int_{\mathbb{R}^{d}}\int_{0}^{T}\Psi_{T}\left(x,s\right)dsdx\right\} \exp\left\{ \int_{\mathbb{R}^{d}}-v_{T}\left(x,0,T\right)dx\right\} ,\label{eq:Laplace_transform_with_w}\end{equation}
 Now we make similar computations to \cite[(3.21)-(3.23)]{BGT2}. By combining
(\ref{eq:Laplace_transform_with_w}) and (\ref{calkowe}) we obtain:\begin{align*}
\mathbb{E}\exp\left\{ -\left\langle \tilde{X}_{T},\Phi\right\rangle \right\}  & =\exp\left\{ \int_{\mathbb{R}^{d}}\int_{0}^{T}\Psi_{T}\left(x,s\right)dsdx\right\} \\
 & \,\,\,\,\,\cdot\exp\left\{ -\int_{\mathbb{R}^{d}}\int_{0}^{T}\Psi_{T}\left(x,T-s\right)\left(1-v_{T}\left(x,T-s,s\right)\right)-VG\left(v_{T}\left(x,T-s,s\right)\right)dsdx\right\} \\
 & =\exp\left\{ \int_{\mathbb{R}^{d}}\int_{0}^{T}\Psi_{T}\left(x,T-s\right)v_{T}\left(x,T-s,s\right)+VG\left(v_{T}\left(x,T-s,s\right)\right)dsdx\right\} \end{align*}
 The last expression can be rewritten as:\begin{equation}
\mathbb{E}\exp\left\{ -\left\langle \tilde{X}_{T},\Phi\right\rangle \right\} =\exp\left\{ V\left(I_{1}\left(T\right)+I_{2}\left(T\right)\right)+I_{3}\left(T\right)\right\} ,\label{eq: Laplace_functional_decomposition}\end{equation}
where \begin{align}
I_{1}\left(T\right)= & \int_{0}^{T}\int_{\mathbb{R}^{d}}\frac{M}{2}\left(\int_{0}^{s}\mathcal{T}_{u}\Psi_{T}\left(\cdot,T+u-s\right)\left(x\right)du\right)^{2}dxds,\nonumber \\
I_{2}\left(T\right)= & \int_{0}^{T}\int_{\mathbb{R}^{d}}\left[G\left(v_{T}\left(x,T-s,s\right)\right)-\frac{M}{2}\left(\int_{0}^{s}\mathcal{T}_{u}\Psi_{T}\left(\cdot,T+u-s\right)\left(x\right)du\right)^{2}\right]dxds,\label{def: I2}\\
I_{3}\left(T\right)= & \int_{0}^{T}\int_{\mathbb{R}^{d}}\Psi_{T}\left(x,T-s\right)v_{T}\left(x,T-s,s\right)dxds.\nonumber \end{align}
To complete the proof we have to compute limits as $T\rightarrow+\infty$.
We claim\begin{equation}
I_{1}\left(T\right)\rightarrow\frac{MK}{2V}\int_{0}^{1}\int_{0}^{1}\int_{\mathbb{R}^{d}}\int_{\mathbb{R}^{d}}\Phi\left(x,t\right)\Phi\left(y,s\right)dxdyC_{h}\left(s,t\right)dsdt,\label{zb:1}\end{equation}
\[
I_{2}\left(T\right)\rightarrow0,\]
\[
I_{3}\left(T\right)\rightarrow0,\]
Combining (\ref{eq: Laplace_functional_decomposition}) with the above
limits we obtain\begin{equation}
\lim_{T\rightarrow+\infty}\mathbb{E}\exp\left\{ -\left\langle \tilde{X}_{T},\Phi\right\rangle \right\} =\exp\left\{ \frac{MK}{2}\int_{0}^{1}\int_{0}^{1}\int_{\mathbb{R}^{d}}\int_{\mathbb{R}^{d}}\Phi\left(x,t\right)\Phi\left(y,s\right)dxdyC_{h}\left(s,t\right)dsdt\right\} \label{eq:limit_<X_T,PHI>}\end{equation}
hence the limit process $X_{T}$ is a Gaussian process with covariance
(\ref{def:kowariancja_tw1}).

\subsubsection{Convergence proofs}

~\\
$I_{1}\left(T\right)$ does not depend on $F$ so it can be evaluated
in the same way as in \cite[(3.32)-(3.34)]{BGT2}. \\
Let us now deal with $I_{3}\left(T\right)$. By using (\ref{ineq: v<n})
we obtain\[
I_{3}\left(T\right)\leq\int_{0}^{T}\int_{\mathbb{R}^{d}}\Psi_{T}\left(x,T-s\right)\int_{0}^{s}\mathcal{T}_{u}\Psi_{T}\left(\cdot,T-u\right)dudxds\leq\]
\[
\frac{C}{F_{T}^{2}}\int_{0}^{T}\int_{\mathbb{R}^{d}}\varphi\left(x\right)\int_{0}^{s}\mathcal{T}_{u}\varphi\left(x\right)dudxds\]
Now the rest of the proof goes along the same lines as in \cite{BGT2}.\\
We will turn to $I_{2}\left(T\right)$ which is a little more intricate.
Combining (\ref{def: I2}) and property $4$ from Fact \ref{fact: properties of G}\begin{align*}
 & I_{2}\left(T\right)=\int_{0}^{T}\int_{\mathbb{R}^{d}}\left[\frac{M}{2}\left[v_{T}\left(\ldots\right)^{2}-\left(\int_{0}^{s}\mathcal{T}_{u}\Psi_{T}\left(\cdot,T+u-s\right)\left(x\right)du\right)^{2}\right]+g\left(v_{T}\left(\ldots\right)\right)v_{T}\left(\ldots\right)^{2}\right]dxds=\\
 & =\frac{M}{2}I'_{2}\left(T\right)+I''_{2}\left(T\right),\end{align*}
where \[
I'_{2}\left(T\right)=\int_{0}^{T}\int_{\mathbb{R}^{d}}v_{T}\left(x,T-s,s\right)^{2}-\left(\int_{0}^{s}\mathcal{T}_{u}\Psi_{T}\left(\cdot,T+u-s\right)\left(x\right)du\right)^{2}dxds,\]
\begin{equation}
I''_{2}\left(T\right)=\int_{0}^{T}\int_{\mathbb{R}^{d}}g\left(v_{T}\left(x,T-s,s\right)\right)v_{T}\left(x,T-s,s\right)^{2}dxds\label{eq: I2pp}\end{equation}
 By inequality (\ref{ineq: v<n}) we have \[
0\leq-I'_{2}\left(T\right)=\int_{0}^{T}\int_{\mathbb{R}^{d}}\left[\left(n_{T}\left(x,T-s,s\right)\right)^{2}-\left(v_{T}\left(x,T-s,s\right)\right)^{2}\right].\]
Combining (\ref{calkowe}) and (\ref{def: n}) yields\begin{align*}
0\leq & n_{T}\left(x,T-s,s\right)-v_{T}\left(x,T-s,s\right)=\\
 & \int_{0}^{s}\mathcal{T}_{s-u}\left[\Psi_{T}\left(\cdot,T-u\right)v_{T}\left(\cdot,T-u,u\right)+VG\left(v_{T}\left(\cdot,T-u,u\right)\right)\right]\left(x\right)du=\left(*\right).\end{align*}

We have $\mathcal{T}_{s}\Psi\geq0$ for $\Psi\geq0$ which is a direct
consequence of the fact that $\mathcal{T}$ is the semigroup of a
Markov process. 
By Fact \ref{G(v)/v2} we have $c(F)$ such that $F(v)\leq \frac{c(F)}{2}v^2$. Hence
\begin{align*}
\left(*\right) & \leq\int_{0}^{s}\mathcal{T}_{s-u}\left[\Psi_{T}\left(\cdot,T-u\right)v_{T}\left(\cdot,T-u,u\right)+c\left(F\right)\frac{V}{2}v_{T}\left(\cdot,T-u,u\right)^{2}\right]\left(x\right)du\\
 & \leq\max\left(1,c\left(F\right)\right)\int_{0}^{s}\mathcal{T}_{s-u}\left[\Psi_{T}\left(\cdot,T-u\right)v_{T}\left(\cdot,T-u,u\right)+\frac{V}{2}v_{T}\left(\cdot,T-u,u\right)^{2}\right]\left(x\right)du\\
 & \leq\max\left(1,c\left(F\right)\right)\int_{0}^{s}\mathcal{T}_{s-u}\left[\Psi_{T}\left(\cdot,T-u\right)n_{T}\left(\cdot,T-u,u\right)+\frac{V}{2}n_{T}\left(\cdot,T-u,u\right)^{2}\right]\left(x\right)du.\end{align*}
Except of the constant $c\left(F\right)$ the last expression does
not depend on $F$.\\
Next we consider\[
n_{T}\left(x,T-s,s\right)+v_{T}\left(x,T-s,s\right)\leq2n_{T}\left(x,T-s,s\right)\leq2\int_{0}^{s}\mathcal{T}_{s-u}\Psi\left(\cdot,T-u\right)\left(x\right)du.\]
The rest of the proof goes along the lines of the proof in \cite[inequalities (3.39)-(3.42)]{BGT2}
and hence we acquire $I_{2}'\left(T\right)\rightarrow0$.\\
Before proving the convergence of $I_{2}''\left(T\right)$ we state
two facts:

\begin{fact}
$n_{T}\left(x,T-s,s\right)\rightarrow0$ in uniformly $x\in\mathbb{R}^{d}$,
$s\in\left[0,T\right]$ as $T\rightarrow+\infty$. \label{fact: uniformaly}
\end{fact}
\begin{proof}
~\\
\[
n_{T}\left(x,T-s,s\right)=\int_{0}^{s}\mathcal{T}_{s-u}\Psi_{T}\left(\cdot,T-u\right)du=\]
\[
\frac{1}{F_{T}}\int_{0}^{s}\mathcal{T}_{s-u}\varphi\left(x\right)\chi\left(\frac{T-u}{T}\right)du\leq\]
\[
\frac{C}{F_{T}}\int_{0}^{+\infty}\mathcal{T}_{u}\varphi\left(x\right)du=\frac{C_{1}}{F_{T}}\int_{\mathbb{R}^{d}}\frac{\varphi\left(y\right)}{\left|x-y\right|^{d-\alpha}}dy\leq\frac{C_{2}}{F_{T}}\rightarrow0.\]
The last line contains the definition of the potential operator of
the semigroup $\mathcal{T}_{t}$ which is bounded in respect to $x$
(this can be found in \cite[Lemma 5.3]{I}).
\end{proof}
\begin{fact}
The following convergence holds:

\[
\int_{0}^{T}\int_{\mathbb{R}^{d}}v_{T}\left(x,T-s,s\right)^{2}\rightarrow c'\left(\Psi\right)\textrm{ as }T\rightarrow+\infty.\]

\end{fact}
\begin{proof}
One easily checks that\[
2\frac{I_{1}\left(T\right)}{M}+I_{2}'\left(T\right)=\int_{0}^{T}\int_{\mathbb{R}^{d}}v_{T}\left(x,T-s,s\right)^{2}.\]
 Hence the result follows from (\ref{zb:1}) and $I'_{2}\left(T\right)\rightarrow0$
as $T\rightarrow0$.
\end{proof}
It is now easy to prove the convergence of $I''_{2}$. From Fact \ref{fact: properties of G}
property $4$ we know that for given $\epsilon>0$ we can choose such
$\delta$ that $\forall_{x\in\left(-\delta,\delta\right)}\left|g\left(x\right)\right|\leq\epsilon$.
Fact \ref{fact: uniformaly} provides us with $T_{0}$ such that $\forall_{T\geq T_{0}}\, n_{T}\left(x,T-s,s\right)<\delta$.
Combining this with (\ref{ineq: v<n}) we obtain $\forall_{T\geq T_{0}}\, g\left(v_{T}\left(x,T-s,s\right)\right)\leq\epsilon$.
Hence for $T>T_{0}$ holds:\begin{align*}
\left|I_{2}''\left(T\right)\right| & \leq\epsilon\int_{0}^{T}\int_{\mathbb{R}^{d}}v_{T}^{2}\left(x,T-s,s\right)dxds\rightarrow\epsilon c'\left(\Psi\right).\end{align*}
Since $\epsilon$ was chosen arbitrary we have convergence: $I_{2}''\left(T\right)\rightarrow0$
hence also $I_{2}\left(T\right)\rightarrow0$ as $T\rightarrow+\infty$.

Thus we obtained the limits for $I_{1},I_{2}$ and $I_{3}$ and the
proof of Theorem \ref{thm:For-the-branching} is completed.

\subsection{Proof of Theorem \ref{thm:equlibiurim}}

\subsubsection{Tightness}

~\\
We begin by claiming that the family $\left\{ X_{T}\right\} _{T>0}$
is tight. Close examination of \cite[Section 3.1]{BGT2} reveals that
only the covariance function of the underlying system is significant
for the proof. By \cite[(3.16)]{BGT1} we know that the covariance function
of the branching system is of the same form as the covariance function of
the non-branching system with the Poisson initial condition. From
this we conclude that $X_{T}$ is tight.

\subsubsection{Laplace functional for $\tilde{X}_{T}$}

~\\
We consider $\tilde{X}_{T}$ defined by (\ref{eq:space-time}). Using (\ref{eq occupation}) and interchanging the order of
integration we obtain\[
\left\langle \tilde{X}_{T},\Phi\right\rangle =\frac{T}{F_{T}}\left[\int_{0}^{1}\left\langle N_{Ts},\Psi\left(\cdot,s\right)\right\rangle ds-\left\langle \lambda,\int_{0}^{1}\Psi\left(\cdot,s\right)ds\right\rangle \right].\]
To prove the convergence of $\tilde{X}_{T}$ to $\tilde{X}$ we will
use its Laplace functional\begin{eqnarray}
 & \mathbb{E}\exp\left\{ -\left\langle \tilde{X}_{T},\Phi\right\rangle \right\} = & \exp\left\{ \int_{\mathbb{R}^{d}}\int_{0}^{T}\Psi_{T}\left(x,t\right)dtdx\right\} \label{eq: Laplace tildeX}\\
 &  & \mathbb{E}\exp\left\{ -\int_{0}^{T}\left\langle N_{s},\Psi_{T}\left(\cdot,s\right)\right\rangle ds\right\} ,\nonumber \end{eqnarray}
 It is easy to check that\begin{equation}
\mathbb{E}\left(\exp\left\{ -\int_{0}^{T}\left\langle N_{s},\Psi_{T}\left(\cdot,s\right)\right\rangle ds\right\} |N_{0}=\mu\right)=\exp\left\{ \left\langle \mu,\ln w_{T}\right\rangle \right\} ,\label{eq: conditionallaplace}\end{equation}
where \[
w_{T}\left(x\right)=\mathbb{E}\exp\left\{ -\int_{0}^{T}\left\langle N_{s}^{x},\Psi_{T}\left(\cdot,t\right)\right\rangle dt\right\} \]
 \\
Now we check that $0\leq-\ln(w_T)$ is integrable. For $T$ big enough by Fact \ref{fact: uniformaly} and inequity (\ref{ineq: v<n}) we have $0\leq v_T \leq c <1$. Hence there exists a constant $C$ such that we have $-\ln(w_T)=-\ln(1-v_T)\leq C v_T \leq C n_T$. A trivial verifications shows that $n_T \in \mathcal{L}^1(\mathbb{R}^d)$ so by (\ref{eq: laplace_equilibrum}) and (\ref{eq: conditionallaplace})
we obtain\[
\mathbb{E}\exp\left\{ -\int_{0}^{T}\left\langle N_{s},\Psi_{T}\left(\cdot,s\right)\right\rangle ds\right\} =\mathbb{E}\left(\mathbb{E}\left(\exp\left\{ -\int_{0}^{T}\left\langle N_{s},\Psi_{T}\left(\cdot,s\right)\right\rangle ds\right\} |N_{0}\right)\right)=\]
\[
\exp\left\{ \left\langle \lambda,w_{T}-1\right\rangle +V\int_{0}^{+\infty}\left\langle \lambda,H\left(W_{T}\left(\cdot,s\right)\right)\right\rangle ds\right\} ,\]
where $W_{T}$ satisfies the equation \[
W_{T}\left(x,l\right)=\mathcal{T}_{l}w_{T}\left(x\right)+V\int_{0}^{l}\mathcal{T}_{l-s}H\left(W_{T}\left(\cdot,s\right)\right)\left(x\right)ds\]
\\
It will be a bit easier to deal with $V_{T}\left(x,l\right)=1-W_{T}\left(x,l\right)$.
The equations have the form (let us recall that $G$ is defined by (\ref{def: G}))
\begin{equation}
	\mathbb{E}\exp\left\{ -\int_{0}^{T}\left\langle N_s,\Psi_{T}\left(\cdot,s\right)\right\rangle ds\right\} =\exp\left\{ \left\langle \lambda,-v_{T}\right\rangle +V\int_{0}^{+\infty}\left\langle \lambda,G\left(V_{T}\left(\cdot,s\right)\right)\right\rangle ds\right\} ,
	\label{eq: inner}
\end{equation}
and \begin{equation}
V_{T}\left(x,l\right)=\mathcal{T}_{l}v_{T}\left(x\right)-V\int_{0}^{l}\mathcal{T}_{l-s}G\left(V_{T}\left(\cdot,s\right)\right)\left(x\right)ds,\label{main_equation}\end{equation}
$W_T$ is defined by (\ref{def: h}) with $\varphi(x)=-\ln w_T(x)$ ($w_T\in[0,1]$ hence $\varphi$ is positive). One can easily see that the definition implies that $W_T\in[0,1]$. Consequently $V_T\in[0,1]$ which together with Fact \ref{G>=3D0} yields $G(V_T)\geq 0$. Hence we obtain an inequality
\begin{equation}
	V_{T}\left(x,l\right)\leq\mathcal{T}_{l}v_{T}\left(x\right),\,\forall_{x\in\mathbb{R}^{d},l\geq0}.
	\label{neq: WT<TvT}
\end{equation}
Combining (\ref{eq: Laplace tildeX}) and (\ref{eq: inner}) we obtain\begin{eqnarray*}
\mathbb{E}\exp\left\{ -\left\langle \tilde{X}_{T},\Phi\right\rangle \right\}  & = & \exp\left\{ \int_{\mathbb{R}^{d}}\int_{0}^{T}\Psi_{T}\left(x,t\right)dtdx\right\} \exp\left\{ -\int_{\mathbb{R}^{d}}v_{T}\left(x\right)dx\right\} \\
 &  & \exp\left\{ V\int_{0}^{+\infty}\int_{\mathbb{R}^{d}}G\left(V_{T}\left(x,t\right)\right)dxdt\right\} =A\left(T\right)\cdot B\left(T\right),\end{eqnarray*}
where \[
A\left(T\right)=\exp\left\{ \int_{\mathbb{R}^{d}}\int_{0}^{T}\Psi_{T}\left(x,t\right)dtdx\right\} \exp\left\{ -\int_{\mathbb{R}^{d}}v_{T}\left(x\right)dx\right\} ,\]
\[
B\left(T\right)=\exp\left\{ V\int_{0}^{+\infty}\int_{\mathbb{R}^{d}}G\left(V_{T}\left(x,t\right)\right)dxdt\right\} .\]
Let us note that $A$ is the same as (\ref{eq:Laplace_transform_with_w})
in the first proof hence we know that its limit is given by (\ref{eq:limit_<X_T,PHI>}).

\subsubsection{Limit of B}

~\\
 To complete the proof the limit $\lim_{T\rightarrow+\infty}B\left(T\right)$
has to be calculated. It suffices to consider 
\begin{equation}
\int_{0}^{+\infty}\int_{\mathbb{R}^{d}}G\left(V_{T}\left(x,t\right)\right)dxdt.\label{exp(...)}\end{equation}
Using Fact \ref{fact: properties of G}, property 4, we split it in
the following way \[
\int_{0}^{+\infty}\int_{\mathbb{R}^{d}}G\left(V_{T}\left(\cdot,t\right)\right)dxdt=\frac{M}{2}\left(B_{1}\left(T\right)+B_{2}\left(T\right)+B_{3}\left(T\right)\right)+B_{4}\left(T\right),\]
 where \[
B_{1}\left(T\right)=\int_{0}^{+\infty}\int_{\mathbb{R}^{d}}V_{T}\left(x,t\right)^{2}-\left(\mathcal{T}_{t}v_{T}\left(x\right)\right)^{2}dxdt,\]
\[
B_{2}\left(T\right)=\int_{0}^{+\infty}\int_{\mathbb{R}^{d}}\left(\mathcal{T}_{t}v_{T}\left(x\right)\right)^{2}-\left(\mathcal{T}_{t}n_{T}\left(x\right)\right)^{2}dxdt,\]
\[
B_{3}\left(T\right)=\int_{0}^{+\infty}\int_{\mathbb{R}^{d}}\left(\mathcal{T}_{t}n_{T}\left(x\right)\right)^{2}dxdt,\]
\[
B_{4}\left(T\right)=\int_{0}^{+\infty}\int_{\mathbb{R}^{d}}g\left(V_{T}\left(x,t\right)\right)V_{T}\left(x,t\right)^{2}.\]
We will prove the following limits (let us recall that we assume (\ref{eq:uproszczenie}) for simplicity)
\[
B_{1}\left(T\right)\rightarrow0,\]
\[
B_{2}\left(T\right)\rightarrow0,\]
\[
B_{3}\left(T\right)\rightarrow\frac{K}{2V}\left\langle \lambda,\varphi\right\rangle ^{2}\int_{0}^{1}\int_{0}^{1}\left\{ -u_{1}^{h}-u_{2}^{h}+\left(u_{1}+u_{2}\right)^{h}\right\} \psi\left(u_{1}\right)\psi\left(u_{2}\right)du_{1}du_{2},\]
\[
B_{4}\left(T\right)\rightarrow0,\]
as $T\rightarrow+\infty$.

\paragraph*{Limit of $B_{1}$}

~\\
By (\ref{neq: WT<TvT}) we obtain \[
0\leq-B_{1}\left(T\right)=\int_{0}^{+\infty}\left(\int_{\mathbb{R}^{d}}\left(\mathcal{T}_{t}v_{T}\left(x\right)\right)^{2}-V_{T}\left(x,t\right)^{2}dx\right)dt=\]
\[
\int_{0}^{+\infty}\int_{\mathbb{R}^{d}}\left(\mathcal{T}_{t}v_{T}\left(x\right)-V_{T}\left(x,t\right)\right)\left(\mathcal{T}_{t}v_{T}\left(x\right)V+_{T}\left(x,t\right)\right)dxdt\leq\]
Combining this with inequality (\ref{neq: WT<TvT}) and equation (\ref{main_equation})
we have\[
\int_{0}^{+\infty}\int_{\mathbb{R}^{d}}\left(V\int_{0}^{t}\mathcal{T}_{t-t'}G\left(V_{T}\left(\cdot,t'\right)\right)\left(x\right)dt'\right)\left(2\mathcal{T}_{t}v_{T}\left(x\right)\right)dxdt=\]
Taking into account the form of $G$ (Fact \ref{fact: properties of G},
property 4)\[
B_{11}\left(T\right)+B_{12}\left(T\right),\]
where \[
B_{11}\left(T\right)=\int_{0}^{+\infty}\int_{\mathbb{R}^{d}}\left(V\frac{M}{2}\int_{0}^{t}\mathcal{T}_{t-t'}V_{T}\left(\cdot,t'\right)^{2}\left(x\right)dt'\right)\left(2\mathcal{T}_{t}v_{T}\left(x\right)\right)dxdt,\]
\[
B_{12}\left(T\right)=\int_{0}^{+\infty}\int_{\mathbb{R}^{d}}\left(V\int_{0}^{t}\mathcal{T}_{t-t'}g\left(V_{T}\left(\cdot,t'\right)\right)V_{T}\left(\cdot,t'\right)^{2}\left(x\right)dt'\right)\left(2\mathcal{T}_{t}v_{T}\left(x\right)\right)dxdt\]
Once again we use inequality (\ref{neq: WT<TvT})\[
B_{11}\left(T\right)\leq VM\int_{0}^{+\infty}\int_{\mathbb{R}^{d}}\left(\int_{0}^{t}\mathcal{T}_{t-t'}\left(\mathcal{T}_{t'}v_{T}\left(\cdot\right)\right)^{2}\left(x\right)dt'\right)\left(\mathcal{T}_{t}v_{T}\left(x\right)\right)dxdt=\]
\[
VM\int_{0}^{+\infty}\int_{0}^{t}\int_{\mathbb{R}^{d}}\mathcal{T}_{t-t'}\left(\mathcal{T}_{t'}v_{T}\left(\cdot\right)\right)^{2}\left(x\right)\mathcal{T}_{t}v_{T}\left(x\right)dxdt'dt\leq\]
Applying (\ref{ineq: v<n}) twice\[
VM\int_{0}^{+\infty}\int_{0}^{t}\int_{\mathbb{R}^{d}}\mathcal{T}_{t-t'}\left(\mathcal{T}_{t'}n_{T}\left(\cdot\right)\right)^{2}\left(x\right)\mathcal{T}_{t}n_{T}\left(x\right)dxdt'dt=\]

\[
MV\int_{0}^{+\infty}\int_{0}^{t}\int_{\mathbb{R}^{d}}\mathcal{T}_{t'}n_{T}\left(x\right)\mathcal{T}_{t'}n_{T}\left(x\right)\mathcal{T}_{2t-t'}n_{T}\left(x\right)dxdt'dt=\]
We use the Plancherel formula and (\ref{eq: Fourier T_t})\[
\frac{MV}{\left(2\pi\right)^{2d}}\int_{0}^{+\infty}\int_{0}^{t}\int_{\mathbb{R}^{2d}}\widehat{\mathcal{T}_{t'}n_{T}}\left(z_{1}\right)\widehat{\mathcal{T}_{t'}n_{T}}(z_{2})\overline{\widehat{\mathcal{T}_{2t-t'}n_{T}}}\left(z_{1}+z_{2}\right)dz_{1}dz_{2}dt'dt=\]
\[
\frac{MV}{\left(2\pi\right)^{2d}}\int_{0}^{+\infty}\int_{0}^{t}\int_{\mathbb{R}^{2d}}e^{-t'\left|z_{1}\right|^{\alpha}}\widehat{n}_{T}\left(z_{1}\right)e^{-t'\left|z_{2}\right|^{\alpha}}\widehat{n}_{T}(z_{2})e^{-(2t-t')\left|z_{1}+z_{2}\right|^{\alpha}}\overline{\widehat{n}_{T}}\left(z_{1}+z_{2}\right)dz_{1}dz_{2}dt'dt=\]
\[
\frac{MV}{\left(2\pi\right)^{2d}}\int_{\mathbb{R}^{2d}}\widehat{n}_{T}\left(z_{1}\right)\widehat{n}_{T}(z_{2})\overline{\widehat{n}_{T}}\left(z_{1}+z_{2}\right)\int_{0}^{+\infty}\int_{0}^{t}e^{-t'\left|z_{1}\right|^{\alpha}}e^{-t'\left|z_{2}\right|^{\alpha}}e^{-(2t-t')\left|z_{1}+z_{2}\right|^{\alpha}}dt'dtdz_{1}dz_{2}=\]
\[
\frac{MV}{\left(2\pi\right)^{2d}}\int_{\mathbb{R}^{2d}}\frac{1}{2\left|z_{1}+z_{2}\right|^{\alpha}\left(\left|z_{1}\right|^{\alpha}+\left|z_{2}\right|^{\alpha}+\left|z_{1}+z_{2}\right|^{\alpha}\right)}\widehat{n}_{T}\left(z_{1}\right)\widehat{n}_{T}(z_{2})\overline{\widehat{n}_{T}}\left(z_{1}+z_{2}\right)dz_{1}dz_{2}=\left(*\right)\]
Before proceeding further we will estimate $\widehat{n}_{T}$
\[
	\left|\widehat{n}_{T}\left(z,r,t\right)\right|=\left|\widehat{\int_{0}^{t}\mathcal{T}_{t-s}\Psi_{T}\left(\cdot,r+t-s\right)ds(z)}\right|=
\]
\[
\left|\frac{1}{F_{T}}\int_{0}^{t}e^{-\left(t-s\right)\left|z\right|^{\alpha}}\widehat{\varphi}\left(z\right)\chi_{T}\left(r+t-s\right)ds\right|\leq\]
\[
\frac{\sup\chi}{F_{T}}\left|\widehat{\varphi}\left(z\right)\right|\int_{0}^{t}e^{-\left(t-s\right)\left|z\right|^{\alpha}}ds\leq\]
 Hence\begin{equation}
\left|\widehat{n}_{T}\left(z,r,t\right)\right|\leq\frac{C}{F_{T}}\frac{\left|\widehat{\varphi}\left(z\right)\right|}{\left|z\right|^{\alpha}}\left[1-e^{-t\left|z\right|^{\alpha}}\right]\label{neq: n(z,r,t)< strong}\end{equation}
and this immediately implies (see (\ref{eq: n_T_x notacja})) \begin{equation}
\left|\widehat{n}_{T}\left(z\right)\right|\leq\frac{C}{F_{T}}\frac{1}{\left|z\right|^{\alpha}}\left[1-e^{-T\left|z\right|^{\alpha}}\right].\label{ineq: hat(n_T)}\end{equation}
Here, and in what follows, $C$ denotes a generic constant. \\
Coming back to $\left(*\right)$ and using the last inequality we
obtain\begin{eqnarray*}
\left|\left(*\right)\right|\leq\frac{C}{F_{T}^{3}}\int_{\mathbb{R}^{2d}}\frac{1}{2\left|z_{1}+z_{2}\right|^{\alpha}\left(\left|z_{1}\right|^{\alpha}+\left|z_{2}\right|^{\alpha}+\left|z_{1}+z_{2}\right|^{\alpha}\right)}\frac{1}{\left|z_{1}\right|^{\alpha}}\left[1-e^{-T\left|z_{1}\right|^{\alpha}}\right]\\
\frac{1}{\left|z_{2}\right|^{\alpha}}\left[1-e^{-T\left|z_{2}\right|^{\alpha}}\right]\frac{1}{\left|z_{1}+z_{2}\right|^{\alpha}}\left[1-e^{-T\left|z_{1}+z_{2}\right|^{\alpha}}\right]dz_{1}dz_{2} & =\end{eqnarray*}
Substituting $T^{1/\alpha}z_{1}=y_{1}$
and $T^{1/\alpha}z_{2}=y_{2}$ yields\begin{eqnarray*}
\frac{CT^{5}}{F_{T}^{3}T^{2\frac{d}{\alpha}}}\int_{\mathbb{R}^{2d}}\frac{1}{\left|y_{1}+y_{2}\right|^{\alpha}\left(\left|y_{1}\right|^{\alpha}+\left|y_{2}\right|^{\alpha}+\left|y_{1}+y_{2}\right|^{\alpha}\right)}\frac{1}{\left|y_{1}\right|^{\alpha}}\left[1-e^{-\left|y_{1}\right|^{\alpha}}\right]\\
\frac{1}{\left|y_{2}\right|^{\alpha}}\left[1-e^{-\left|y_{2}\right|^{\alpha}}\right]\frac{1}{\left|y_{1}+y_{2}\right|^{\alpha}}\left[1-e^{-\left|y_{1}+y_{2}\right|^{\alpha}}\right]dy_{1}dy_{2} & \leq\end{eqnarray*}
\[
B'_{11}\left(T\right)\cdot B''_{11},\]
where\[
B'_{11}\left(T\right)=\frac{C'T^{5}}{F_{T}^{3}T^{2\frac{d}{\alpha}}}\]
\begin{eqnarray*}
B''_{11}=\int_{\mathbb{R}^{2d}}\frac{1}{\left|y_{1}+y_{2}\right|^{\alpha}\left(\left|y_{1}\right|^{\alpha}+\left|y_{2}\right|^{\alpha}+\left|y_{1}+y_{2}\right|^{\alpha}\right)}\frac{1}{\left|y_{1}\right|^{\alpha}}\left[1-e^{-\left|y_{1}\right|^{\alpha}}\right]\\
\frac{1}{\left|y_{2}\right|^{\alpha}}\left[1-e^{-\left|y_{2}\right|^{\alpha}}\right]\frac{1}{\left|y_{1}+y_{2}\right|^{\alpha}}\left[1-e^{-\left|y_{1}+y_{2}\right|^{\alpha}}\right]dy_{1}dy_{2}\end{eqnarray*}
The integral $B''_{11}$ is finite which will be proved in Fact \ref{fact: calka1}.
The expression $B'_{11}\left(T\right)$ can be evaluated\[
B'_{11}\left(T\right)=T^{\frac{10-3\left(3-\frac{d}{\alpha}\right)-4\frac{d}{\alpha}}{2}}=T^{\frac{1-\frac{d}{\alpha}}{2}}\]
and as $1-\frac{d}{\alpha}<0$ $B'_{11}\left(T\right)\rightarrow0$
hence the convergence: $B_{11}\left(T\right)\rightarrow0$ is obtained
too.\\
From Fact \ref{fact: uniformaly} and inequalities (\ref{ineq: v<n})
and (\ref{neq: WT<TvT}) we know $V_{T}(x,l)\rightarrow0$
uniformly as $T\rightarrow0$ and so $g\left(V_{T}\left(x,l\right)\right)\leq\epsilon$
for $T$ sufficiently large hence \[
B_{12}\left(T\right)\leq\epsilon\int_{0}^{+\infty}\int_{\mathbb{R}^{d}}\left(V\int_{0}^{t}\mathcal{T}_{t-t'}V_{T}\left(\cdot,t'\right)^{2}\left(x\right)dt'\right)\left(2\mathcal{T}_{t}v_{T}\left(x\right)\right)dxdt\leq\frac{2\epsilon}{M}B_{11}\left(T\right)\]
 thus $B_{12}\left(T\right)\rightarrow0$ and $B_{1}\left(T\right)\rightarrow0$
too.

\paragraph*{Limit of $B_{2}$}

~\\
Let us first estimate expression $n_{T}-v_{T}$ using (\ref{calkowe})
and (\ref{def: n}) \begin{eqnarray*}
 & n_{T}\left(x\right)-v_{T}\left(x\right)= & \int_{0}^{T}\mathcal{T}_{T-u}\Psi_{T}\left(\cdot,T-u\right)\left(x\right)du\\
 &  & -\int_{0}^{T}\mathcal{T}_{T-u}\left[\Psi_{T}\left(\cdot,T-u\right)\left(1-v_{T}\left(\cdot,T-u,u\right)\right)-VG\left(v_{T}\left(\cdot,T-u,u\right)\right)\right]\left(x\right)du\end{eqnarray*}
\[
n_{T}\left(x\right)-v_{T}\left(x\right)=\int_{0}^{T}\mathcal{T}_{T-u}\left[\Psi_{T}\left(\cdot,T-u\right)v_{T}\left(\cdot,T-u,u\right)+VG\left(v_{T}\left(\cdot,T-u,u\right)\right)\right]\left(x\right)du\leq\]
Applying Fact \ref{G(v)/v2}\[
\int_{0}^{T}\mathcal{T}_{T-u}\left[\Psi_{T}\left(\cdot,T-u\right)v_{T}\left(\cdot,T-u,u\right)+Vc\left(v_{T}\left(\cdot,T-u,u\right)\right)^{2}\right]\left(x\right)du,\]
where $c$ is a constant. By inequality (\ref{ineq: v<n})\begin{equation}
n_{T}\left(x\right)-v_{T}\left(x\right)\leq\int_{0}^{T}\mathcal{T}_{T-u}\left[\Psi_{T}\left(\cdot,T-u\right)n_{T}\left(\cdot,T-u,u\right)+Vc\left(n_{T}\left(\cdot,T-u,u\right)\right)^{2}\right]\left(x\right)du\label{ineq: n-v}\end{equation}
We have\[
0\leq-B_{2}\left(T\right)=\int_{0}^{+\infty}\left(\int_{\mathbb{R}^{d}}\left(\mathcal{T}_{t}n_{T}\left(x\right)\right)^{2}-\left(\mathcal{T}_{t}v_{T}\left(x\right)\right)^{2}dx\right)dt=\]
\[
\int_{0}^{+\infty}\int_{\mathbb{R}^{d}}\left(\mathcal{T}_{t}\left(n_{T}\left(\cdot\right)-v_{T}\left(\cdot\right)\right)\left(x\right)\right)\left(\mathcal{T}_{t}\left(v_{T}\left(\cdot\right)+n_{T}\left(\cdot\right)\right)\left(x\right)\right)dxdt\leq\]
Applying (\ref{ineq: v<n}) and (\ref{ineq: n-v})\begin{eqnarray*}
2\int_{0}^{+\infty}\int_{\mathbb{R}^{d}}\mathcal{T}_{t}\left\{ \int_{0}^{T}\mathcal{T}_{T-u}\left[\Psi_{T}\left(\cdot,T-u\right)n_{T}\left(\cdot,T-u,u\right)+Vc\left(n_{T}\left(\cdot,T-u,u\right)\right)^{2}\right]du\right\} \left(x\right)\\
\mathcal{T}_{t}n_{T}\left(x\right)dxdt=\end{eqnarray*}
Now we apply the Plancherel formula\begin{eqnarray*}
\frac{2}{\left(2\pi\right)^{d}}\int_{0}^{+\infty}\int_{\mathbb{R}^{d}}e^{-2t\left|z\right|^{\alpha}}\widehat{\int_{0}^{T}\mathcal{T}_{T-u}\left[\Psi_{T}\left(\cdot,T-u\right)n_{T}\left(\cdot,T-u,u\right)+Vc\left(n_{T}\left(\cdot,T-u,u\right)\right)^{2}\left(\cdot\right)\right]\left(z\right)du}\\
\widehat{n}_{T}\left(z\right)dzdt & =\end{eqnarray*}
Interchanging the order of integration and integrating with respect
to $t$ we get\begin{eqnarray*}
\frac{1}{\left(2\pi\right)^{d}}\int_{\mathbb{R}^{d}}\frac{1}{\left|z\right|^{\alpha}}\widehat{\int_{0}^{T}\mathcal{T}_{T-u}\left[\Psi_{T}\left(\cdot,T-u\right)n_{T}\left(\cdot,T-u,u\right)+Vc\left(n_{T}\left(\cdot,T-u,u\right)\right)^{2}\left(\cdot\right)\right]\left(z\right)du}\\
\widehat{n}_{T}\left(z\right)dz=c'\left(B_{21}\left(T\right)+B_{22}\left(T\right)\right),\end{eqnarray*}
where \[
B_{21}\left(T\right)=\int_{\mathbb{R}^{d}}\frac{1}{\left|z\right|^{\alpha}}\left\{ \widehat{\int_{0}^{T}\mathcal{T}_{T-u}\left[\Psi_{T}\left(\cdot,T-u\right)n_{T}\left(\cdot,T-u,u\right)\right]\left(z\right)du}\right\} \widehat{n}_{T}\left(z\right)dz,\]
\[
B_{22}\left(T\right)=\int_{\mathbb{R}^{d}}\frac{1}{\left|z\right|^{\alpha}}\left\{ \widehat{\int_{0}^{T}\mathcal{T}_{T-u}\left[Vc\left(n_{T}\left(\cdot,T-u,u\right)\right)^{2}\left(\cdot\right)\right]\left(z\right)du}\right\} \widehat{n}_{T}\left(z\right)dz.\]
We shall compute $\lim_{T\rightarrow+\infty}B_{21}\left(T\right)$
first. We have
\[
B_{21}(T)=\int_{\mathbb{R}^{d}}\frac{1}{|z|^{\alpha}} \left\{ \int_{0}^{T} e^{-(T-u)|z|^\alpha} 
\widehat{\Psi_{T}\left(\cdot,T-u\right)}*\widehat{n_{T}\left(\cdot,T-u,u\right)\left(z\right)}
\right\} \widehat{n}_T(z)dz.
\]
The inner convolution can be estimated using inequality (\ref{neq: n(z,r,t)< strong}) and simplification (\ref{eq:uproszczenie})
\[
\left|\widehat{\Psi_{T}\left(\cdot,T-u\right)}*\widehat{n_{T}\left(\cdot,T-u,u\right)}\left(z\right)\right|=
\left|\chi_{T}\left(T-u\right)\widehat{\varphi_{T}}\left(\cdot\right)*\widehat{n_{T}\left(\cdot,T-u,u\right)\left(z\right)}\right|=\]
\[
\left|\chi_{T}\left(T-u\right)\int_{\mathbb{R}^{d}}\widehat{\varphi}_{T}\left(z-x\right)\widehat{n}_{T}\left(x,T-u,u\right)dx\right|\leq\]
\[
\frac{c\left(\chi\right)}{F_{T}^{2}}\chi_{T}\left(T-u\right)\int_{\mathbb{R}^{d}}\left|\widehat{\varphi}\left(z-x\right)\widehat{\varphi}\left(x\right)\right|\frac{1}{\left|x\right|^{\alpha}}dx\leq\frac{C}{F_{T}^{2}}\]
In the last inequality we use the fact that $\hat{\varphi}$ is bounded
and $\frac{\widehat{\varphi}\left(x\right)}{\left|x\right|^{\alpha}}$
is integrable. Hence we have inequality\begin{equation}
\left|\widehat{\Psi_{T}\left(\cdot,T-u\right)}*\widehat{n_{T}\left(\cdot,T-u,u\right)}\left(z\right)\right|\leq\frac{C}{F_{T}^{2}}\label{neq:convultion1}\end{equation}
 Thus $B_{21}$ satisfies\[
\left|B_{21}\left(T\right)\right|\leq\frac{C}{F_{T}^{2}}\int_{\mathbb{R}^{d}}\frac{1}{\left|z\right|^{\alpha}}\int_{0}^{T}e^{-\left(T-u\right)\left|z\right|^{\alpha}}du\cdot\widehat{n}_{T}\left(z\right)dz\leq\]
Using inequality (\ref{ineq: hat(n_T)}) and integrating with respect
to $u$\[
C'\frac{1}{F_{T}^{3}}\int_{\mathbb{R}^{d}}\frac{1}{\left|z\right|^{\alpha}}\frac{1}{\left|z\right|^{\alpha}}\left[1-e^{-T\left|z\right|^{\alpha}}\right]\frac{1}{\left|z\right|^{\alpha}}\left[1-e^{-T\left|z\right|^{\alpha}}\right]dz=\]
Substituting $zT^{1/\alpha}=y$\[
C'\frac{T^{3}}{F_{T}^{3}T^{\frac{d}{\alpha}}}\int_{\mathbb{R}^{d}}\frac{1}{\left|y\right|^{\alpha}}\frac{1}{\left|y\right|^{\alpha}}\left[1-e^{-\left|y\right|^{\alpha}}\right]\frac{1}{\left|y\right|^{\alpha}}\left[1-e^{-\left|y\right|^{\alpha}}\right]dy\leq B_{21}'\left(T\right)\cdot B''_{21},\]
where\[
B_{21}'\left(T\right)=C''\frac{T^{3}}{F_{T}^{3}T^{\frac{d}{\alpha}}}\]
\[
B_{21}''=\int_{\mathbb{R}^{d}}\frac{1}{\left|y\right|^{\alpha}}\frac{1}{\left|y\right|^{\alpha}}\left[1-e^{-\left|y\right|^{\alpha}}\right]\frac{1}{\left|y\right|^{\alpha}}\left[1-e^{-\left|y\right|^{\alpha}}\right]dy\]
Is is clear that integral $B_{21}''$ in the
last expression is finite since in a neighborhood of
$0$ the integrated expression is proportional to $\frac{1}{\left|y\right|^{\alpha}}$
and it is $O\left(\frac{1}{\left|y\right|^{3\alpha}}\right)$ as $\left|y\right|\rightarrow+\infty$ 
(recall that $\alpha<d<2\alpha$). Now only $B_{21}'$ needs to be evaluated\[
B_{21}'\left(T\right)=C''T^{\frac{6-3\left(3-\frac{d}{\alpha}\right)-2\frac{d}{\alpha}}{2}}=C''T^{\frac{-3+\frac{d}{\alpha}}{2}}.\]
Hence it is obvious that $B_{21}'\left(T\right)\rightarrow0$ as $T\rightarrow0$
and so $\lim_{T\rightarrow0}B_{21}\left(T\right)=0$.\\
Before proceeding to $B_{22}$ we will make the following estimation
using inequality (\ref{neq: n(z,r,t)< strong}) \[
\left|\widehat{\left(n_{T}\left(\cdot,T-u,u\right)\right)^{2}}\right|\left(z\right)=\left|\int_{\mathbb{R}^{d}}\widehat{n}_{T}\left(x,T-u,u\right)\widehat{n}_{T}\left(z-x,T-u,u\right)dx\right|\leq\]
\[
\frac{C}{F_{T}^{2}}\int_{\mathbb{R}^{d}}\frac{1}{\left|x\right|^{\alpha}}\left[1-e^{-u\left|x\right|^{\alpha}}\right]\frac{1}{\left|z-x\right|^{\alpha}}\left[1-e^{-u\left|z-x\right|^{\alpha}}\right]dx\leq\]
 Substitution $xu^{1/\alpha}=y$ yields\[
u^{2-\frac{d}{\alpha}}\frac{C}{F_{T}^{2}}\int_{\mathbb{R}^{d}}\frac{1}{\left|y\right|^{\alpha}}\left[1-e^{-\left|y\right|^{\alpha}}\right]\frac{1}{\left|zu^{1/\alpha}-y\right|^{\alpha}}\left[1-e^{-\left|zu^{1/\alpha}-y\right|^{\alpha}}\right]dy\leq\frac{C'}{F_{T}^{2}}u^{2-\frac{d}{\alpha}}\]
since the integral can be regarded as a convolution of $\mathcal{L}^{2}$
functions so it is bounded. This clearly implies\[
\left|B_{22}\left(T\right)\right|\leq\frac{C'}{F_{T}^{2}}\int_{\mathbb{R}^{d}}\frac{1}{\left|z\right|^{\alpha}}\int_{0}^{T}e^{-\left(T-u\right)\left|z\right|^{\alpha}}u^{2-\frac{d}{\alpha}}du\cdot\left|\widehat{n}_{T}\left(z\right)\right|dz\leq\]
\[
C'\frac{T^{2-\frac{d}{\alpha}}}{F_{T}^{2}}\int_{\mathbb{R}^{d}}\frac{1}{\left|z\right|^{\alpha}}\int_{0}^{T}e^{-\left(T-u\right)\left|z\right|^{\alpha}}du\cdot\left|\widehat{n}_{T}\left(z\right)\right|dz\leq\]
Using inequality (\ref{ineq: hat(n_T)}) we obtain\[
C''\frac{T^{2-\frac{d}{\alpha}}}{F_{T}^{3}}\int_{\mathbb{R}^{d}}\frac{1}{\left|z\right|^{\alpha}}\frac{1}{\left|z\right|^{\alpha}}\left(1-e^{-T\left|z\right|^{\alpha}}\right)\frac{1}{\left|z\right|^{\alpha}}\left(1-e^{-T\left|z\right|^{\alpha}}\right)dz=\]
Substituting $zT^{1/\alpha}=y$ we can rewrite the last expression
as\[
C''\frac{T^{5-\frac{d}{\alpha}}}{F_{T}^{3}T^{\frac{d}{\alpha}}}\int_{\mathbb{R}^{d}}\frac{1}{\left|y\right|^{\alpha}}\frac{1}{\left|y\right|^{\alpha}}\left(1-e^{-\left|y\right|^{\alpha}}\right)\frac{1}{\left|y\right|^{\alpha}}\left(1-e^{-\left|y\right|^{\alpha}}\right)dy.\]
The integral is finite (the same proof as for $B_{21}''$) and \[
\frac{T^{5-\frac{d}{\alpha}}}{F_{T}^{3}T^{\frac{d}{\alpha}}}=T^{\frac{10-2\frac{d}{\alpha}-3\left(3-\frac{d}{\alpha}\right)-2\frac{d}{\alpha}}{2}}=T^{\frac{1-\frac{d}{\alpha}}{2}}\]
 which yields $B_{22}\left(T\right)\rightarrow0$ as $T\rightarrow+\infty$

\paragraph*{Limit of $B_{3}$}

~\\
Applying the Plancherel formula to $B_{3}\left(T\right)$ we get\[
B_{3}\left(T\right)=\frac{1}{\left(2\pi\right)^{d}}\int_{0}^{\infty}\int_{\mathbb{R}^{d}}e^{-2t\left|z\right|^{\alpha}}\left(\widehat{n}_{T}\left(z\right)\right)^{2}dzdt=\]
\[
\frac{1}{\left(2\pi\right)^{d}}\int_{\mathbb{R}^{d}}\left(\widehat{n}_{T}\left(z\right)\right)^{2}\int_{0}^{\infty}e^{-2t\left|z\right|^{\alpha}}dtdz=\]
\[
\frac{1}{2\left(2\pi\right)^{d}}\int_{\mathbb{R}^{d}}\frac{1}{\left|z\right|^{\alpha}}\left(\widehat{n}_{T}\left(z\right)\right)^{2}dz=\]
\[
\frac{1}{2\left(2\pi\right)^{d}}\int_{\mathbb{R}^{d}}\frac{1}{\left|z\right|^{\alpha}}\left(\int_{0}^{T}e^{-\left(T-u\right)\left|z\right|^{\alpha}}\widehat{\varphi}_{T}\left(z\right)\chi_{T}\left(T-u\right)du\right)^{2}dz=\]
\[
\frac{1}{2\left(2\pi\right)^{d}}\frac{1}{F_{T}^{2}}\int_{\mathbb{R}^{d}}\frac{1}{\left|z\right|^{\alpha}}\left(\int_{0}^{T}e^{-u\left|z\right|^{\alpha}}\widehat{\varphi}\left(z\right)\chi_{T}\left(u\right)du\right)^{2}dz=\]
Substituting $u'=u/T$\[
\frac{1}{2\left(2\pi\right)^{d}}\frac{T^{2}}{F_{T}^{2}}\int_{\mathbb{R}^{d}}\frac{1}{\left|z\right|^{\alpha}}\left(\int_{0}^{1}e^{-Tu'\left|z\right|^{\alpha}}\widehat{\varphi}\left(z\right)\chi\left(u'\right)du'\right)^{2}dz=\]
\[
\frac{1}{2\left(2\pi\right)^{d}}\frac{T^{2}}{F_{T}^{2}}\int_{0}^{1}\int_{0}^{1}\int_{\mathbb{R}^{d}}\frac{1}{\left|z\right|^{\alpha}}e^{-T\left(u_{1}+u_{2}\right)\left|z\right|^{\alpha}}\left(\widehat{\varphi}\left(z\right)\right)^{2}\chi\left(u_{1}\right)\chi\left(u_{2}\right)du_{1}du_{2}dz=\]
Let $z=\left[T\left(u_{1}+u_{2}\right)\right]^{-\frac{1}{\alpha}}y$\begin{eqnarray*}
\frac{1}{2\left(2\pi\right)^{d}}\frac{T^{3-\frac{d}{\alpha}}}{F_{T}^{2}}\int_{0}^{1}\int_{0}^{1}\int_{\mathbb{R}^{d}}\left(u_{1}+u_{2}\right)\frac{1}{\left|y\right|^{\alpha}}e^{-\left|y\right|^{\alpha}}\left(\widehat{\varphi}\left(\left[T\left(u_{1}+u_{2}\right)\right]^{-\frac{1}{\alpha}}y\right)\right)^{2}\\
\left(u_{1}+u_{2}\right)^{-\frac{d}{\alpha}}\chi\left(u_{1}\right)\chi\left(u_{2}\right)du_{1}du_{2}dy\end{eqnarray*}
Therefore by Lebesgue's dominated convergence theorem we obtain the
limit of $B_{3}\left(T\right)$\[
\lim_{T\rightarrow+\infty}B_{3}\left(T\right)=\frac{1}{2\left(2\pi\right)^{d}}\int_{0}^{1}\int_{0}^{1}\int_{\mathbb{R}^{d}}\left(u_{1}+u_{2}\right)^{1-\frac{d}{\alpha}}\frac{1}{\left|y\right|^{\alpha}}e^{-\left|y\right|^{\alpha}}\left(\widehat{\varphi}\left(0\right)\right)^{2}\chi\left(u_{1}\right)\chi\left(u_{2}\right)du_{1}du_{2}dy=\]
\[
\frac{\Gamma\left(\frac{d}{\alpha}-1\right)}{2^{d}\alpha\Gamma\left(\frac{d}{2}\right)\pi^{\frac{d}{2}}}\left\langle \lambda,\varphi\right\rangle ^{2}\int_{0}^{1}\int_{0}^{1}\left(u_{1}+u_{2}\right)^{1-\frac{d}{\alpha}}\chi\left(u_{1}\right)\chi\left(u_{2}\right)du_{1}du_{2}=\]
Integrating by parts\[
\frac{K}{2V}\left\langle \lambda,\varphi\right\rangle ^{2}\int_{0}^{1}\int_{0}^{1}\left\{ -u_{1}^{h}-u_{2}^{h}+\left(u_{1}+u_{2}\right)^{h}\right\} \psi\left(u_{1}\right)\psi\left(u_{2}\right)du_{1}du_{2}\]

\paragraph*{Limit of $B_{4}$}

~\\
Firstly, let us notice that \[
B_{1}\left(T\right)+B_{2}\left(T\right)+B_{3}\left(T\right)=\int_{0}^{+\infty}\int_{\mathbb{R}^{d}}V_{T}\left(x,t\right)^{2},\]
 and hence \[
\int_{0}^{+\infty}\int_{\mathbb{R}^{d}}V_{T}\left(x,t\right)^{2}\rightarrow_{T\rightarrow+\infty}C.\]
Secondly by Fact \ref{fact: uniformaly} and inequalities (\ref{neq: WT<TvT})
and (\ref{ineq: v<n}) we know $V_{T}\left(x\right)\rightarrow0$
uniformly as $T\rightarrow0$. Hence $g\left(W_{T}\left(x\right)\right)\leq\epsilon$
for $T$ sufficiently large so \[
\left|B_{4}\left(T\right)\right|\leq\epsilon\int_{0}^{+\infty}\int_{\mathbb{R}^{d}}V_{T}\left(x,t\right)^{2},\]
which clearly implies that $B_{4}\left(T\right)\rightarrow0$ as $T\rightarrow+\infty$.

\paragraph*{Putting the results together}

~\\
Combining the previous results we conclude\[
\lim_{T\rightarrow+\infty}B\left(T\right)=\exp\left\{ \frac{MK}{4}\left\langle \lambda,\varphi\right\rangle ^{2}\int_{0}^{1}\int_{0}^{1}\left\{ -u_{1}^{h}-u_{2}^{h}+\left(u_{1}+u_{2}\right)^{h}\right\} \psi\left(u_{1}\right)\psi\left(u_{2}\right)du_{1}du_{2}\right\} \]
And finally by (\ref{eq:limit_<X_T,PHI>})\[
\lim_{T\rightarrow+\infty}A\left(T\right)B\left(T\right)=\exp\left\{ \frac{MK}{2}\left\langle \lambda,\varphi\right\rangle ^{2}\int_{0}^{1}\int_{0}^{1}c_{h}\left(u_{1},u_{2}\right)\psi\left(u_{1}\right)\psi\left(u_{2}\right)du_{1}du_{2}\right\} ,\]
where $c_{h}$ is the covariance function of fractional Brownian
motion defined by (\ref{def: c - fractional}). This Laplace functional
defines a process $\tilde{X}_{T}$ corresponding to the Gaussian process
$X_{T}$ with the covariance (\ref{def:kowariancja_tw2}) hence Theorem
\ref{thm:equlibiurim} is proved.

\section{Appendix}

The appendix contains a technical fact used in the main proof.

\begin{fact}
\label{fact: calka1} \begin{eqnarray*}
 &  & \int_{\mathbb{R}^{2d}}\frac{1}{\left|y_{1}+y_{2}\right|^{\alpha}\left(\left|y_{1}\right|^{\alpha}+\left|y_{2}\right|^{\alpha}+\left|y_{1}+y_{2}\right|^{\alpha}\right)}\frac{1}{\left|y_{1}\right|^{\alpha}}\left[1-e^{-\left|y_{1}\right|^{\alpha}}\right]\\
 &  & \frac{1}{\left|y_{2}\right|^{\alpha}}\left[1-e^{-\left|y_{2}\right|^{\alpha}}\right]\frac{1}{\left|y_{1}+y_{2}\right|^{\alpha}}\left[1-e^{-\left|y_{1}+y_{2}\right|^{\alpha}}\right]dy_{1}dy_{2}<+\infty\end{eqnarray*}

\end{fact}
\begin{proof}
Substituting $x=y_{1}+y_{2}$ and $z=y_{2}$ we get\[
\int_{\mathbb{R}^{2d}}\frac{1}{\left|x\right|^{\alpha}\left(\left|x\right|^{\alpha}+\left|z\right|^{\alpha}+\left|x-z\right|^{\alpha}\right)}\frac{1}{\left|x-z\right|^{\alpha}}\left[1-e^{-\left|x-z\right|^{\alpha}}\right]\frac{1}{\left|z\right|^{\alpha}}\left[1-e^{-\left|z\right|^{\alpha}}\right]\frac{1}{\left|x\right|^{\alpha}}\left[1-e^{-\left|x\right|^{\alpha}}\right]dxdz=\]
\[
\int_{\mathbb{R}^{2d}}\frac{1}{\left|x\right|^{\alpha}}\frac{1}{\left|x\right|^{\alpha}}\left[1-e^{-\left|x\right|^{\alpha}}\right]\int_{\mathbb{R}^{d}}\frac{1}{\left|x\right|^{\alpha}+\left|z\right|^{\alpha}+\left|x-z\right|^{\alpha}}\frac{1}{\left|x-z\right|^{\alpha}}\left[1-e^{-\left|x-z\right|^{\alpha}}\right]\frac{1}{\left|z\right|^{\alpha}}\left[1-e^{-\left|z\right|^{\alpha}}\right]dzdx=\left(*\right)\]
Let us investigate now
\[
\int_{\mathbb{R}^{d}}\frac{1}{\left|x\right|^{\alpha}+\left|z\right|^{\alpha}+\left|x-z\right|^{\alpha}}\frac{1}{\left|x-z\right|^{\alpha}}\left[1-e^{-\left|x-z\right|^{\alpha}}\right]\frac{1}{\left|z\right|^{\alpha}}\left[1-e^{-\left|z\right|^{\alpha}}\right]dz \leq%\label{def:int_in_fact1}%
\]
\[
 \int_{\mathbb{R}^{d}}\frac{1}{\left|z\right|^{\alpha}}\frac{1}{\left|z\right|^{\alpha}}\left[1-e^{-\left|z\right|^{\alpha}}\right]\frac{1}{\left|x-z\right|^{\alpha}}\left[1-e^{-\left|x-z\right|^{\alpha}}\right]dz \leq c\int_{\mathbb{R}^{d}}\frac{1}{\left|z\right|^{\alpha}}\frac{1}{\left|z\right|^{\alpha}}\left[1-e^{-\left|z\right|^{\alpha}}\right]dz
\]

% For positive $a,b$ the following inequality holds \[
% \frac{1}{a+b}\leq\frac{1}{a}+\frac{1}{b}.\]
% Hence\[
% \frac{1}{\left|x\right|^{\alpha}+\left|z\right|^{\alpha}+\left|x-z\right|^{\alpha}}\leq\frac{1}{\left|z\right|^{\alpha}+\left|x-z\right|^{\alpha}}\leq\frac{1}{\left|z\right|^{\alpha}}+\frac{1}{\left|x-z\right|^{\alpha}},\]
% so (\ref{def:int_in_fact1}) can be estimated\begin{eqnarray*}
%  & \leq & \int_{\mathbb{R}^{d}}\frac{1}{\left|z\right|^{\alpha}}\frac{1}{\left|z\right|^{\alpha}}\left[1-e^{-\left|z\right|^{\alpha}}\right]\frac{1}{\left|x-z\right|^{\alpha}}\left[1-e^{-\left|x-z\right|^{\alpha}}\right]dz+\\
%  &  & \int_{\mathbb{R}^{d}}\frac{1}{\left|x-z\right|^{\alpha}}\frac{1}{\left|x-z\right|^{\alpha}}\left[1-e^{-\left|x-z\right|^{\alpha}}\right]\frac{1}{\left|z\right|^{\alpha}}\left[1-e^{-\left|z\right|^{\alpha}}\right]dz=\end{eqnarray*}
% \[
% 2\int_{\mathbb{R}^{d}}\frac{1}{\left|z\right|^{\alpha}}\frac{1}{\left|z\right|^{\alpha}}\left[1-e^{-\left|z\right|^{\alpha}}\right]\frac{1}{\left|x-z\right|^{\alpha}}\left[1-e^{-\left|x-z\right|^{\alpha}}\right]dz\leq c\int_{\mathbb{R}^{d}}\frac{1}{\left|z\right|^{\alpha}}\frac{1}{\left|z\right|^{\alpha}}\left[1-e^{-\left|z\right|^{\alpha}}\right]dz=c_{2}\]
The last integral is finite since in the neighborhood of $0$ the
integrated function is $O\left(\frac{1}{\left|z\right|^{\alpha}}\right)$
and for big $\left|z\right|$ is $O\left(\frac{1}{\left|z\right|^{2\alpha}}\right)$.
Going back to $\left(*\right)$ we obtain\[
\left(*\right)\leq c_{2}\int_{\mathbb{R}^{d}}\frac{1}{\left|x\right|^{\alpha}}\frac{1}{\left|x\right|^{\alpha}}\left[1-e^{-\left|x\right|^{\alpha}}\right]<c_{3},\]
by the same reason as above.
\end{proof}
\begin{acknowledgement*}
The author would like to thank his supervisor - prof. Tomasz Bojdecki
- for much appreciated help given in general introduction to the branching
systems theory and in writing this paper. The author wishes to thank also prof. Luis Gorostiza for several helpful comments.
\end{acknowledgement*}


\begin{thebibliography}{}
\bibitem[1]{B}P. Billingsley, Convergence of Probability Measures., John Wiley\&Sons, New York, 1968.
\bibitem[2]{BGR}T. Bojdecki, L.G. Gorostiza and S. Ramaswami, %1986,%
 Convergence of $\mathcal{S}'$-valued processes and space time random fields,
J. Funct. Anal. 66 (1986), pp. 21-41.
\bibitem[3]{BGT1}T. Bojdecki, L.G. Gorostiza and A. Talarczyk, Sub-fractional Brownian
motion and its relation to occupation times, Statist. Probab. Lett. 69 (2004), pp. 405-419.
\bibitem[4]{BGT2}T. Bojdecki, L.G. Gorostiza and A. Talarczyk, Limit theorems for
occupation time fluctuations of branching systems I: Long-range dependence, Stoch. Proc. Appl. 116 (2006), pp. 1-18.
\bibitem[5]{BGT3}T. Bojdecki, L.G. Gorostiza and A. Talarczyk, 
Limit theorems for occupation time fluctuations of branching systems II: Critical and large dimensions
 Functional, Stoch. Proc. Appl. 116 (2006), pp. 19-35.
\bibitem[6]{BGT4}T. Bojdecki, L.G. Gorostiza and A. Talarczyk, %2005a.% 
A long range dependence stable process and an infinite variance branching system, 
www.arxiv.org, math.PR/0511739 (2005).
\bibitem[7]{BGT5}T. Bojdecki, L.G. Gorostiza and A. Talarczyk, %2005b.% 
Occupation time fluctuations of an infinite variance branching systems in large
dimensions, www.arxiv.org, math.PR/0511745 (2005).
\bibitem[8]{BZ}M. Birkner and I. Z\"ahle, %2005,% 
Functional central limit theorems for the occupation time of the origin for branching random walks in
$d\geq3$, Weierstra\ss{} Insitut f\"ur Angewandte Analysis und Stochastik,
Berlin, preprint No. 1011 (2005).
\bibitem[9]{DW}J.D. Deuschel and K. Wang, %(1994)%, 
Large deviations for the occupation time of a Poisson system of independent Brownian particles, Stoch. Proc. Appl. 52 (1994), pp. 183-209.
\bibitem[10]{GR}L.G. Gorositza and E.R. Rodrigues, %1999% 
A stochastic model for transport of particulate matter in air: an asymptotic analysis. Acta
Appl. Math. 59 (1999), pp. 21-43.
\bibitem[11]{GW1}L.G. Gorostiza and A. Wakolbinger, Long time behavior of critical branching 
particle systems and its applications, CRM Proc. and Lect. Notes Vol. 5 (1994), pp. 119-137.
\bibitem[12]{GW2}L.G. Gorostiza and A. Wakolbinger, %1991.% 
Persistence criteria for a class of critical branching particle systems in continuous time.
Ann. Probab. 19 (1991), pp. 266-288.
\bibitem[13]{I}I. Iscoe, A weighted occupation time for a class of measure-valued
branching processes, Probab. Th. Rel. Fields 71 (1986), pp. 85-116.
\bibitem[14]{M}I. Mitoma, Tightness of probabilities on $C\left(\left[0,1\right],\mathcal{S}'\right)$
and $D\left(\left[0,1\right],\mathcal{S}'\right)$, Ann. Probab. 11 (1983), pp. 989-999.
\end{thebibliography}
\end{document}